\documentstyle{article}

\newcommand{\qed}{\hfill $\Box$ \\[0.2cm]}

\def\Tn{{\mbox{${\cal T}_n$}}}

\def\Jn{{\mbox{${\cal J}_n$}}}

\newcommand{\Z}{{\mbox{\bf Z}}}

\newcommand{\N}{{\mbox{\bf N}}}

\newcommand{\mj}{\mbox{\bf 1}}

\begin{document}

\title
{\bf A New Proof of the Faithfulness of Brauer's
Representation of Temperley-Lieb Algebras}

\author{
{\sc Kosta Do\v sen}$ ^\ast$, {\sc \v Zana Kovijani\' c}$ ^{\ast\ast}$
and {\sc Zoran Petri\' c}$ ^\ast$
\\[.3cm]
$ ^\ast$Matemati\v cki institut, SANU
\\Knez Mihailova 35, p.f. 367
\\11001 Belgrade, Serbia and Montenegro
\\email: \{kosta, zpetric\}@mi.sanu.ac.yu
\\[.2cm]
$ ^{\ast\ast}$Prirodno-matemati\v cki fakultet
\\Univerzitet Crne Gore
\\Cetinjski put b.b., p.f. 211
\\81000 Podgorica, Serbia and Montenegro
\\email: zanak@cg.yu
}

\date{ }
\maketitle
\begin{abstract}
\noindent The faithfulness of the orthogonal group case of
Brauer's representation of the Brauer centralizer algebras
restricted to their Temperley-Lieb subalgebras, which was
established by Vaughan Jones, is here proved in a new, elementary
and self-contained, manner.

\vspace{0.3cm}

\\[.2cm]
\\[.2cm]
{\it Acknowledgement}. We would like to thank an anonymous referee
for a very helpful comment concerning the beginning of Section~4.
\end{abstract}

\section{Introduction}

Temperley-Lieb algebras, whose roots are in
statistical mechanics (see \cite{TL71} and \cite{GHJ89}, Appendix II.b),
play a prominent role in knot theory and low-dimensional topology. They
have entered this field through Jones' representation of Artin's braid
groups in a version of these algebras (see \cite{J87}, \S11, and works
of Jones cited therein), on which the famous Jones polynomial of
knot theory is based (see also \cite{KL}, \cite{L97}, \cite{PS}). The
version of Temperley-Lieb algebra we deal with here, which is pretty
standard (see \cite{KL}, and works by Kauffman and others cited therein),
is defined as follows.

The Temperley-Lieb algebra \Tn, where $n\geq 2$, is an associative algebra
over a field
of scalars whose nature is not important for us in this paper (in
\cite{KL}, p. 9, it is taken to be the field of rational functions
$P/Q$ with $P,Q\in\Z[x,x^{-1}]$). The basis of \Tn, which we will call \Jn,
is freely generated from a set of generators $\mj^n,h^n_1,\ldots,h^n_{n-1}$
with multiplication subject to the monoid equations (i.e., multiplication is
associative and $\mj^n$ is a unit for it) and the equations
\[
\begin{array}{ll}
\makebox[3em][l]{(1)} & \makebox[15em][l]{$h^n_ih^n_j=h^n_jh^n_i$,
\quad $i-j\geq 2$,}
\\[0.2cm]
(2) & h^n_ih^n_{i\pm 1}h^n_i=h^n_i,
\\[0.2cm]
(3) & h^n_ih^n_i=h^n_i.
\end{array}
\]
For the definition of multiplication in \Tn\ the equation $(3)$ is replaced
by
\[
\begin{array}{ll}
\makebox[3em][l]{$(3p)$} & \makebox[15em][l]{$h^n_ih^n_i=ph^n_i$}
\end{array}
\]
with $p$ a specified nonzero scalar (in \cite{KL}, it is taken to be
$-x^2-x^{-2}$, but other values are found in other treatments of Temperley-Lieb
algebras). The basis \Jn\ of \Tn\ is finite: its cardinality is the
$n$-th Catalan number $(2n)!/(n!(n+1)!)$. The algebra \Tn\ is the vector space
whose basis is \Jn.

For those \Tn\ where $p$ in $(3p)$ is a natural number
greater than or equal to 2 there is a representation in
matrices due to Brauer (see \cite{B37} and \cite{W88}, Section~3).
This is the orthogonal group case of Brauer's representation
restricted to the Temperley-Lieb subalgebra of the Brauer
algebra. The faithfulness, i.e.\ isomorphism, of this representation
was established by Jones in \cite{J94} (Section~3, Theorem 3.4, p. 330)
by referring to
his technique of Markov trace. Our purpose in this paper is
to give a different, elementary and self-contained, proof of that
faithfulness. We believe our proof is worth publishing because of
its aesthetic value.

The computational interest of Brauer's representation is lessened by the fact that for
\Tn\ we have to pass to $p^n\times p^n$ matrices, and so get an exponential
growth. However, this growth is to be expected, since for sufficiently
large $n$ the $n$-th Catalan number is greater than $(2-\varepsilon)^{2n}$
(this can be computed with the help of Stirling's Theorem). When they were first
introduced in \cite{TL71}, Temperley-Lieb algebras were represented in matrices
in a manner different from Brauer's (see \cite{GHJ89}, Appendix II.b, p. 264),
again with an exponentional growth.

Brauer's representation is provided by assigning to the elements of the basis
\Jn\ of \Tn\ particular 0-1 matrices that satisfy (1), (2) and (3$p$)
for $p$ a natural number greater than or equal 2. For the faithfulness of the
representation it suffices to establish that the list
of matrices assigned to the elements of \Jn\ is linearly independent.
Our proof of linear independence proceeds as follows. We introduce a linear order
on the matrices, and establish that every matrix has an entry with 1 where all
the matrices preceding it in the order have 0.

It is shown in \cite{DP} that Brauer's representation of
Temperley-Lieb algebras is based on
the fact that the Kronecker product of matrices gives rise to
an endofunctor of the category of matrices that is adjoint
to itself. The result of this paper is interesting because of
applications in the area covered by \cite{DP}.

\section{Ordering \Jn}

For $n\geq 2$ and $1\leq i\leq n-1$, let $h^n_{i,1}$ be $h^n_i$, and for
$2\leq k\leq i$ let $h^n_{i,k}$ be $h^n_{i,k-1}h^n_{i-k+1}$. We call the
expressions $h^n_{i,k}$, for $1\leq k\leq i$, {\it blocks}.

Every element of \Jn\ different from $\mj^n$ is denoted by a unique expression in
{\it Jones' normal form}
\[ h^n_{i_1,k_1}\ldots h^n_{i_l,k_l} \]
where $l\geq 1$, $i_1<\ldots< i_l$ and
$i_1-k_1<\ldots<i_l-k_l$ (see \cite{J83}, \S 4.1.4, \cite{GHJ89}, \S 2.8,
\cite{BDP} and \cite{DP}, \S 10). Let $h^n_{i,k}$ stand for $\mj^n$ when $k=0$,
and let us call these expressions too {\it blocks}. Then every element of \Jn\ is
denoted by a unique expression in the {\it normal form}
\[h^n_{1,k_1}\ldots h^n_{n-1,k_{n-1}}\]
where $0\leq k_i\leq i$, and for every $k_i,k_j>0$ such that $i<j$
we have $i-k_i<j-k_j$. So we may identify every element of \Jn\ by a sequence
$k_1\ldots k_{n-1}$
subject to the conditions on $k_i$ and $k_j$ we have just stated.

Let these sequences be ordered lexicographically from the right; in
other words, let
\[k_1'\ldots k_{n-1}'<k_1\ldots k_{n-1}\]
if and only if for some $i$ we have $k_i'< k_i$ and for every $j>i$ we have $k_j'=k_j$.
This induces a linear order on \Jn, which we call the
{\em lexicographical order} of \Jn.

\section{Matrices, relations and graphs for \Jn}

Let $p$ be a natural number greater than or equal to 2, and let $E_p$ be the
$1\times p^2$ matrix that for $1\leq i,j\leq p$ has the entries
\[E_p(1,(i-1)p+j)=\delta(i,j),\]
where $\delta$ is the Kronecker delta. For example, $E_2$ is $[1\;0\;0\;1]$
and $E_3$ is $[1\;0\;0\;0\;1\;0\;0\;0\;1]$. Let $E_p'$ be the transpose of $E_p$,
and let us assign to $h^n_k$, $1\leq k\leq n-1$, the $p^n\times p^n$ matrix
$I_{p^{n-k-1}}\otimes E_p' E_p\otimes I_{p^{k-1}}$, where $I_m$ is the
$m\times m$ identity matrix with the entries $I_m(i,j)=\delta(i,j)$,
and $\otimes$ is the Kronecker product of matrices
(see \cite{J53}, Chapter VII.5, pp. 211-213). We assign to $\mj^n$ the matrix
$I_{p^n}$.

Every $n\times m$ matrix $A$ whose entries are only 0 and 1 may be identified with
a binary relation $R_A\subseteq n\times m$ such that $A(i,j)=1$ if and only if
$(i,j)\in R_A$.
Every binary relation may of course be drawn as a bipartite graph. Here are a
few examples of such graphs for matrices we have introduced up to now, with
$p=2$:

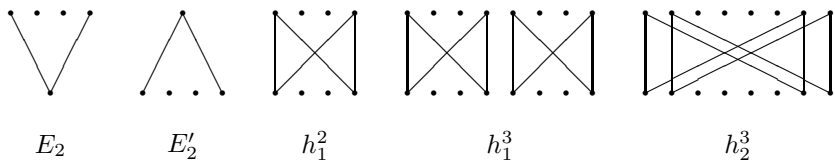
\begin{center}
\begin{picture}(310,90)
\put(0,70){\circle*{2}}
\put(10,70){\circle*{2}}
\put(20,70){\circle*{2}}
\put(30,70){\circle*{2}}
\put(65,70){\circle*{2}}
\put(100,70){\circle*{2}}
\put(110,70){\circle*{2}}
\put(120,70){\circle*{2}}
\put(130,70){\circle*{2}}
\put(150,70){\circle*{2}}
\put(160,70){\circle*{2}}
\put(170,70){\circle*{2}}
\put(180,70){\circle*{2}}
\put(190,70){\circle*{2}}
\put(200,70){\circle*{2}}
\put(210,70){\circle*{2}}
\put(220,70){\circle*{2}}
\put(240,70){\circle*{2}}
\put(250,70){\circle*{2}}
\put(260,70){\circle*{2}}
\put(270,70){\circle*{2}}
\put(280,70){\circle*{2}}
\put(290,70){\circle*{2}}
\put(300,70){\circle*{2}}
\put(310,70){\circle*{2}}
\put(15,40){\circle*{2}}
\put(50,40){\circle*{2}}
\put(60,40){\circle*{2}}
\put(70,40){\circle*{2}}
\put(80,40){\circle*{2}}
\put(100,40){\circle*{2}}
\put(110,40){\circle*{2}}
\put(120,40){\circle*{2}}
\put(130,40){\circle*{2}}
\put(150,40){\circle*{2}}
\put(160,40){\circle*{2}}
\put(170,40){\circle*{2}}
\put(180,40){\circle*{2}}
\put(190,40){\circle*{2}}
\put(200,40){\circle*{2}}
\put(210,40){\circle*{2}}
\put(220,40){\circle*{2}}
\put(240,40){\circle*{2}}
\put(250,40){\circle*{2}}
\put(260,40){\circle*{2}}
\put(270,40){\circle*{2}}
\put(280,40){\circle*{2}}
\put(290,40){\circle*{2}}
\put(300,40){\circle*{2}}
\put(310,40){\circle*{2}}

\put(0,70){\line(1,-2){15}}
\put(30,70){\line(-1,-2){15}}

\put(50,40){\line(1,2){15}}
\put(80,40){\line(-1,2){15}}

\put(100,40){\line(0,1){30}}
\put(100,40){\line(1,1){30}}
\put(130,40){\line(-1,1){30}}
\put(130,40){\line(0,1){30}}

\put(150,40){\line(0,1){30}}
\put(150,40){\line(1,1){30}}
\put(180,40){\line(-1,1){30}}
\put(180,40){\line(0,1){30}}
\put(190,40){\line(0,1){30}}
\put(190,40){\line(1,1){30}}
\put(220,40){\line(-1,1){30}}
\put(220,40){\line(0,1){30}}

\put(240,40){\line(0,1){30}}
\put(240,40){\line(2,1){60}}
\put(300,40){\line(-2,1){60}}
\put(300,40){\line(0,1){30}}
\put(250,40){\line(0,1){30}}
\put(250,40){\line(2,1){60}}
\put(310,40){\line(-2,1){60}}
\put(310,40){\line(0,1){30}}

\put(15,20){\makebox(0,0){$E_2$}}
\put(65,20){\makebox(0,0){$E'_2$}}
\put(115,20){\makebox(0,0){$h^2_1$}}
\put(185,20){\makebox(0,0){$h^3_1$}}
\put(275,20){\makebox(0,0){$h^3_2$}}

\end{picture}
\end{center}

If $\mj^n, h^n_1, \ldots,h^n_{n-1}$ denote the 0-1 matrices
we have assigned to
these expressions, then it can be verified that these matrices
satisfy the equations
(1), (2) and (3$p$), with multiplication being matrix multiplication
(see \cite{DP}, \S\S 17-20). If
$\mj^n,h^n_1,\ldots,h^n_{n-1}$ denote the corresponding binary relations,
then for multiplication being composition of binary relations the
equations (1), (2) and (3) are satisfied. In both cases we also have the monoid
equations.

Composition of binary relations is easy to read
from bipartite graphs. Here is an example:

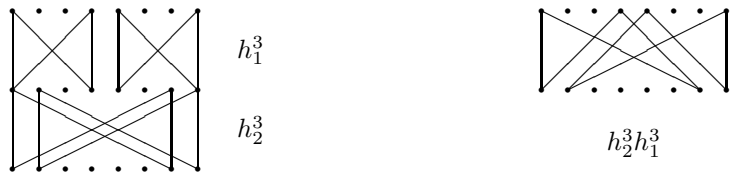
\begin{center}
\begin{picture}(270,100)
\put(0,20){\circle*{2}}
\put(10,20){\circle*{2}}
\put(20,20){\circle*{2}}
\put(30,20){\circle*{2}}
\put(40,20){\circle*{2}}
\put(50,20){\circle*{2}}
\put(60,20){\circle*{2}}
\put(70,20){\circle*{2}}

\put(0,50){\circle*{2}}
\put(10,50){\circle*{2}}
\put(20,50){\circle*{2}}
\put(30,50){\circle*{2}}
\put(40,50){\circle*{2}}
\put(50,50){\circle*{2}}
\put(60,50){\circle*{2}}
\put(70,50){\circle*{2}}

\put(0,80){\circle*{2}}
\put(10,80){\circle*{2}}
\put(20,80){\circle*{2}}
\put(30,80){\circle*{2}}
\put(40,80){\circle*{2}}
\put(50,80){\circle*{2}}
\put(60,80){\circle*{2}}
\put(70,80){\circle*{2}}

\put(200,50){\circle*{2}}
\put(210,50){\circle*{2}}
\put(220,50){\circle*{2}}
\put(230,50){\circle*{2}}
\put(240,50){\circle*{2}}
\put(250,50){\circle*{2}}
\put(260,50){\circle*{2}}
\put(270,50){\circle*{2}}

\put(200,80){\circle*{2}}
\put(210,80){\circle*{2}}
\put(220,80){\circle*{2}}
\put(230,80){\circle*{2}}
\put(240,80){\circle*{2}}
\put(250,80){\circle*{2}}
\put(260,80){\circle*{2}}
\put(270,80){\circle*{2}}

\put(0,50){\line(0,1){30}}
\put(0,50){\line(1,1){30}}
\put(30,50){\line(-1,1){30}}
\put(30,50){\line(0,1){30}}
\put(40,50){\line(0,1){30}}
\put(40,50){\line(1,1){30}}
\put(70,50){\line(-1,1){30}}
\put(70,50){\line(0,1){30}}

\put(0,20){\line(0,1){30}}
\put(0,20){\line(2,1){60}}
\put(60,20){\line(-2,1){60}}
\put(60,20){\line(0,1){30}}
\put(10,20){\line(0,1){30}}
\put(10,20){\line(2,1){60}}
\put(70,20){\line(-2,1){60}}
\put(70,20){\line(0,1){30}}

\put(200,50){\line(0,1){30}}
\put(200,50){\line(1,1){30}}
\put(210,50){\line(1,1){30}}
\put(210,50){\line(2,1){60}}
\put(260,50){\line(-2,1){60}}
\put(260,50){\line(-1,1){30}}
\put(270,50){\line(-1,1){30}}
\put(270,50){\line(0,1){30}}

\put(90,35){\makebox(0,0){$h^3_2$}}
\put(90,65){\makebox(0,0){$h^3_1$}}
\put(235,30){\makebox(0,0){$h^3_2h^3_1$}}

\end{picture}
\end{center}

By so composing binary relations we can assign to every element of \Jn\ a binary
relation, and then from this binary relation we can recover the 0-1 matrix assigned
to our element of \Jn. For this to provide an isomorphic representation of \Tn\
in matrices it is sufficient (and also necessary) to establish that the list of
0-1 matrices assigned to the elements of \Jn\ is linearly independent. The
remainder of this paper is devoted to establishing this fact.

\section{Diagonals}

We may use $\mj^n,h^n_1,\ldots,h^n_{n-1}$, and expressions obtained from these
by multiplying, to denote either the 0-1 matrices
assigned to these expressions in
the previous section, or the corresponding binary relations, or the corresponding
bipartite graphs. So we may speak of the 1 entries of $h^n_i$, conceived as a
matrix, which correspond to the ordered pairs of $h^n_i$, conceived as a binary
relation, which correspond to the edges of $h^n_i$, conceived as a bipartite
graph. When formulating our results we will stick, however, to the
terminology of binary relations (but it helps intuition if examples of such
relations are drawn as bipartite graphs).

For $n, k\geq 0$, let $A^n_k$ be the set
$\{k+1,\ldots,k+p^n\}\subseteq \N$, and for $1\leq q\leq p$ let
\[S_q A^{n+1}_k=\{(q-1)p^n+j \mid j\in A^n_k\}.\]
For example, if $p=2$, then $A^2_0=\{1,2,3,4\}$, $A^1_0=\{1,2\}$,
$A^1_2=\{3,4\}$, $A^0_2=\{3\}$, $S_1A^2_0=A^1_0$, $S_2A^2_0=A^1_2$
and $S_1S_2A^2_0=S_1(S_2A^2_0)=A^0_2$. Words in the alphabet
$\{S_q \mid 1\leq q\leq p\}$ will be called $S$-{\it words}.

For $W$ a word, let $W^0$  be the empty word, and let $W^{n+1}$ be
$W^n W$. Let $|W|$ be the length of the word $W$. We will use the letters
$W, V, U, W_1,\ldots,W',\ldots$ for $S$-words.

For $1\leq q_1, q_2\leq p$, and $q_1\neq q_2$, let $L$ be short for
$S_{q_1}$ and $R$ for $S_{q_2}$. (When $p=2$, in bipartite graphs $L=S_1$
is interpreted as ``left'', and $R=S_2$ as ``right''.) Then
the following is an easy consequence of definitions.
\\[.2cm]
{\sc Remark 1.}\quad {\it If $n\geq 2$ and $|V|=n-2$, then the
pairs $(VLLA^n_0,VRRA^n_0)$ and $(VRRA^n_0,VLLA^n_0)$ are in the
binary relation $h^n_{n-1}$.} \vspace{.2cm}

From now on, we will abbreviate $(W_1A^n_0,W_2A^n_0)$ to
$(W_1,W_2)$, omitting $A^n_0$, and, moreover, we will take for
granted that such pairs belong to a binary relation, without
specifying it explicitly. Then the following remark generalizes
Remark~1.
\\[.2cm]
{\sc Remark 2.}\quad {\it If $|V|\geq 0$, $|W|\geq 0$, $i=1+|V|$ and
$n=i+1+|W|$, then $(VLLW,VRRW)$ and $(VRRW,VLLW)$ are in $h^n_i$.}
\vspace{.2cm}

Next we establish the following lemma (due to
the second author), which is fundamental for our proof.
\\[.2cm]
{\sc Lemma 3.}\quad {\it If $|V|\geq 0$, $|W|\geq 0$, $k\geq 0$,
$i=k+|V|\geq 1$ and $n=i+1+|W|$, then in $h^n_{i,k}$ we have pairs of the
following forms:}
\\[.2cm]
\makebox[3em][l]{({\it even})} {\it for $k=2l$,}
\[
\begin{array}{ll}
\makebox[5em][l]{(I {\it even})} &
\makebox[20em][l]{$(V(RL)^l LW \; , \; VL(LR)^l W)$,}
\\
\makebox[5em][l]{(II {\it even})} &
\makebox[20em][l]{$(V(LR)^l RW \; , \; VR(RL)^l W)$;}
\end{array}
\]
\makebox[3em][l]{({\it odd})} {\it for $k=2l+1$,}
\[
\begin{array}{ll}
\makebox[5em][l]{(I {\it odd})} &
\makebox[20em][l]{$(V(LR)^l LLW \; , \; VRR(LR)^l W)$,}
\\
\makebox[5em][l]{(II {\it odd})} &
\makebox[20em][l]{$(V(RL)^l RRW \; , \; VLL(RL)^l W)$.}
\end{array}
\]
(Note that (II {\it even}) is obtained from (I {\it even}) by interchanging
$L$ and $R$, and
analogously for (II {\it odd}) and (I {\it odd}). Note also that
in ({\it even}) the
$S$-words in between $V$ and $W$ on the right-hand sides are obtained from the
corresponding words on the left-hand sides by reading them in reverse order, while
in ({\it odd}) we have to read them in reverse order and interchange $L$ and $R$.)
\\[.2cm]
{\sc Proof of Lemma 3.}\quad If $k=0$, then in $h^n_{i,0}$, which is $\mj^n$,
we have $(VLW,VLW)$ and $(VRW,VRW)$. For $1\leq k\leq i$ we proceed
by induction on $k$ in $h^n_{i,k}$. In the basis, when $k=1$, we apply
Remark~2. For the induction step, when $2\leq k\leq i$, we have
\[
h^n_{i,k}=h^n_{i,k-1}h^n_{i-k+1},\]
and suppose the Lemma holds for $h^n_{i,k-1}$.

If $k-1=2l+1$, then for
$|V|=i-k$ and $|W|=n-i-1$ from (I {\it odd}) of the induction hypothesis we
obtain that
\[
(VR(LR)^l LLW,VRRR(LR)^l W)\]
is in $h^n_{i,k-1}$, and from
Remark~2 we obtain that
\[
(VRRR(LR)^l W,VLLR(LR)^l W)\]
is in
$h^n_{i-k+1}$, which yields (I {\it even}) for $h^n_{i,k}$.
We obtain analogously
(II {\it even}) for $h^n_{i,k}$ from (II {\it odd}) for $h^n_{i,k-1}$
and Remark~2
(we just interchange $L$ and $R$).

If $k-1=2l+2$, then for $|V|=i-k$ and $|W|=n-i-1$ from (I {\it even}) of the
induction hypothesis we obtain that
\[
(VL(RL)^{l+1} LW,VLL(LR)^{l+1} W)\]
is in $h^n_{i,k-1}$, and from Remark~2
we obtain that
\[
(VLL(LR)^{l+1} W,RR(LR)^{l+1}W)\]
is in $h^n_{i-k+1}$,
which yields (I {\it odd}) for $h^n_{i,k}$. We obtain analogously (II {\it odd})
for $h^n_{i,k}$ from (II {\it even}) for $h^n_{i,k-1}$ and Remark~2.
\qed

The pairs mentioned in Lemma~3 will be called the {\em diagonals} of the block
$h^n_{i,k}$ for $1\leq k\leq i$. As particular cases of these diagonals
we have the pairs mentioned in Remarks 1 and 2. For $k=0$, where
$h^n_{i,0}$ is $\mj^n$, all the pairs of $\mj^n$ are the {\em diagonals}
of $\mj^n$.

Let us say that a nonintersecting {\it n}-diagram $D$ of the kind
introduced by Kauffman, and considered in \cite{BDP} and \cite{DP},
can be {\it put} into the pair $(W_1,W_2)$ when the threads of $D$ join
the same $S$-symbols in the $S$-word $W_1W_2$. Then if $i=k>0$ and
$n=i+1$, one can put into the diagonals of $h^n_{i,k}$ just the
{\it n}-diagram corresponding to this block. An analogous property holds
for the ``diagonals'' of other matrices, which will be defined below.

We can establish the following.
\\[.2cm]
{\sc Remark 4.}\quad {\it If $|V|\geq 0$, $|W|\geq 0$, $k\geq 2$,
$i=k+|V|$, $n=i+1+|W|$, $1\leq k'<k$ and $|U|=n$, then
in $h^n_{i,k'}$ we don't have any pair of the following forms:}
\\[.2cm]
\makebox[3em][l]{({\it odd})} {\it for $k=2l+1$,
$(U,VRR(LR)^l W)$ and $(U,VLL(RL)^l W)$,}
\\[.1cm]
\makebox[3em][l]{({\it even})} {\it for $k=2l+2$,
$(U,VL(LR)^{l+1}W)$ and $(U,VR(RL)^{l+1}W)$.}
\\[.2cm]
This is an easy consequence of the definition of $h^n_{i,k}$. It suffices to look
at the right members of the diagonals of $h^n_{i,k}$, which cannot occur as
right members of any pair of $h^n_{i,k'}$ for $1\leq k'<k$.
\\[.2cm]
{\sc Remark 5.}\quad {\it If $|V|\geq 0$,
$|W|\geq 0$, $k\geq 0$, $i=k+|V|\geq 1$,
$n=i+1+|W|$, $|U|=n$ and $|U'|=i$, then}
\\[.1cm]
\makebox[5em][l]{({\it even})} {\it for} $k=2l$,
\\[.1cm]
\makebox[5em][l]{(I {\it even})} {\it every pair $(U,VL(LR)^{l}W)$ of
$h^n_{i,k}$ is either a diagonal of $h^n_{i,k}$ or

\makebox[3.5em][l]{ } $U$ is of the form
$U'S_qW$ for a symbol $S_q$ different from $L$;}
\\[.05cm]
\makebox[5em][l]{(II {\it even})} {\it every pair $(U,VR(RL)^{l}W)$ of
$h^n_{i,k}$ is either a diagonal of $h^n_{i,k}$ or

\makebox[3.5em][l]{ } $U$ is of the form
$U'S_qW$ for a symbol $S_q$ different from $R$;}
\\[.1cm]
\makebox[5em][l]{({\it odd})} {\it for} $k=2l+1$,
\\[.1cm]
\makebox[5em][l]{(I {\it odd})} {\it every pair $(U,VRR(LR)^l W)$ of
$h^n_{i,k}$ is either a diagonal of $h^n_{i,k}$ or

\makebox[3.5em][l]{ } $U$ is of the form
$U'S_qW$ for a symbol $S_q$ different from $L$;}
\\[.05cm]
\makebox[5em][l]{(II {\it odd})} {\it every pair $(U,VLL(RL)^l W)$ of
$h^n_{i,k}$ is either a diagonal of $h^n_{i,k}$ or

\makebox[3.5em][l]{ } $U$ is of the form
$U'S_qW$ for a symbol $S_q$ different from $R$.}
\\[.2cm]
This too is an easy consequence of the definition of $h^n_{i,k}$.

Let $t$ be in the normal form
\[h^n_{1,k_1}\ldots h^n_{n-1,k_{n-1}}\]
of Section~2, and let $\max(t)=\max\{i\mid k_i>0\}$. If for
every $i$ we have $k_i=0$, then $\max(t)$ is undefined, and $t$ is equal to $\mj^n$.
\\[.2cm]
{\sc Remark 6.}\quad {\em If $\max(t)=i$ and $n=i+1+w$,
then for every $W$ with $|W|=w$ there are $V$ and $U$ with
$|V|=|U|=i+1$ such that in $t$ we have $(VW,UW)$,
and there are no other pairs in $t$ except those of such a form.}
\\[.2cm]
This is an easy consequence of the definition of $h^n_i$.

If $t$ is in normal form as above, then we define the {\it diagonals} of $t$
inductively in terms of the diagonals of blocks. Formally, we define the
diagonals of $h^n_{1,k_1}\ldots h^n_{i,k_i}$ for every $i$ such that
$1\leq i\leq n-1$. The diagonals of $h^n_{1,k_1}$ are the diagonals of this
block. If $(W_1,W_2)$ is a diagonal of $h^n_{1,k_1}\ldots h^n_{j,k_j}$
with $j<n-1$ and $(W_2,W_3)$ is a diagonal of $h^n_{j+1,k_{j+1}}$, then
$(W_1,W_3)$
is a diagonal of $h^n_{1,k_1}\ldots h^n_{j+1,k_{j+1}}$.
\\[.2cm]
{\sc Remark 7.}\quad {\it If $\max(t)=i$ and $n=i+1+|W|$ for $|W|\geq 0$,
then every diagonal of $t$ is of the form $(VLW,URW)$ or of the form
$(VRW,ULW)$ for some $V$ and $U$ such that
$|V|=|U|=i$.}
\vspace{.2cm}

From the definition above it is not clear that every $t$ has diagonals.
If $t$ is equal to $\mj^n$, then it certainly has diagonals. For the rest we have the
following lemma.
\\[.2cm]
{\sc Lemma 8.}\quad {\it Consider $h^n_{1,k_1}\ldots h^n_{n-1,k_{n-1}}$ in
normal form. If $t$ is $h^n_{1,k_1}\ldots h^n_{j,k_j}$ for $j\leq n-1$,
$|V|\geq 0$, $|W|\geq 0$, $\max(t)=i=|V|+k_i\geq 1$,
$n=i+1+|W|$ and $|U|=n$,
then, for some $V$ and every $W$, in $t$ we have diagonals of the following
forms:}
\\[.2cm]
\makebox[3em][l]{({\it even})} {\it for $k_i=2l$,}
\[
\begin{array}{ll}
\makebox[5em][l]{(I {\it even})} &
\makebox[20em][l]{$(U \; , \; VL(LR)^l W)$,}
\\
\makebox[5em][l]{(II {\it even})} &
\makebox[20em][l]{$(U \; , \; VR(RL)^l W)$;}
\end{array}
\]
\makebox[3em][l]{({\it odd})} {\it for $k_i=2l+1$,}
\[
\begin{array}{ll}
\makebox[5em][l]{(I {\it odd})} &
\makebox[20em][l]{$(U \; , \; VRR(LR)^l W)$,}
\\
\makebox[5em][l]{(II {\it odd})} &
\makebox[20em][l]{$(U \; , \; VLL(RL)^l W)$.}
\end{array}
\]
{\sc Proof.}\quad We proceed by induction on $i$. For $i=1$ we
apply Lemma~3. For $i\geq 2$ let $t$ be $t'h^n_{i,k_i}\ldots h^n_{j,k_j}$,
and suppose the lemma holds for $t'$. If $t'$ is equal to $\mj^n$, it
suffices to apply Lemma~3 to $h^n_{i,k_i}$. Otherwise, let $\max(t')=i'$.
We have to consider four cases.
\\[.1cm]
({\it even}-{\it even})\quad If $k_{i'}=2l_{i'}$ and $k_i=2l_i$, then,
by the induction
hypothesis, for $|V'|=i'-k_{i'}$ and $|W'|=n-i'-1$ in $t'$ we have
\[
\begin{array}{ll}
(\mbox{\rm I}'\; \mbox{\it even}) & (U,V'L(LR)^{l_{i'}}W'),
\\
(\mbox{\rm II}'\; \mbox{\it even}) & (U,V'R(RL)^{l_{i'}}W'),
\end{array}
\]
and by Lemma~3 for $|V|=i-k_i$ and $|W|=n-i-1$ in $h^n_{i,k_i}$
we have
\[
\begin{array}{ll}
(\mbox{\rm I}\; \mbox{\it even}) & (V(RL)^{l_i}LW,VL(LR)^{l_{i}}W),
\\
(\mbox{\rm II}\; \mbox{\it even}) & (V(LR)^{l_i}RW,VR(RL)^{l_{i}}W).
\end{array}
\]

Since $t$ is in normal form, we have $i'<i$ and $i'-k_{i'}<i-k_i$,
and so we have $|W|<|W'|$ and
$|V'|<|V|$. Then we can replace $W'$ in
$(\mbox{\rm I}'\; \mbox{\it even})$ and $(\mbox{\rm II}'\; \mbox{\it even})$
by $W''LW$ and $W''RW$ for $|W''|\geq 0$, and we can replace $V$
in (I {\it even}) and (II {\it even}) by $V'LV''$ and $V'RV''$ for $|V''|\geq 0$.
We have $|V''|=|W''|+2(l_{i'}-l_i)$, and so $|V''|$ is even if and only if
$|W''|$ is even. If both $|V''|$ and $|W''|$ are even,
then we make $(LR)^{l_{i'}}W''$ equal to $V''(LR)^{l_i}$, and
$(RL)^{l_{i'}}W''$ equal to $V''(RL)^{l_i}$. If both $|V''|$ and $|W''|$
are odd, then we make $(LR)^{l_{i'}}W''$ equal to $V''(RL)^{l_i}$, and
$(RL)^{l_{i'}}W''$ equal to $V''(LR)^{l_i}$.

We proceed analogously in the remaining three cases, where $k_{i'}$
is even and $k_i$ odd, where $k_{i'}$ is odd and $k_i$ even,
and, finally, where they are both odd.
\qed

\section{Linear independence for \Jn}

We are now ready to prove the following result, which guarantees linear
independence for the representation of \Jn.
\\[.2cm]
{\sc Theorem.}\quad {\it If $t'$ precedes $t$ in the lexicographical order,
then the diagonals of $t$ are not pairs of $t'$.}
\\[.2cm]
{\sc Proof.}\quad Let $t$ be $h^n_{1,k_1}\ldots h^n_{n-1,k_{n-1}}$, let
$t'$ be $h^n_{1,k'_1}\ldots h^n_{n-1,k'_{n-1}}$, and let $t'$ precede $t$
in the lexicographical order. Then there is an $i$ such that $k'_i<k_i$
and for every $j>i$ we have $k'_j=k_j$. We make an induction on the number
$\nu$ of $j$'s such that $j>i$ and $k'_j=k_j>0$.

If $\nu$ is 0, then $\max(t)=i$. If $t'$ is $\mj^n$, then the Theorem follows
by Remark~7. If $t'$ is not $\mj^n$, then let $\max(t')=i'$. Hence, $k_i$ and
$k'_i$ are both greater than 0, and we must have $i'\leq i$. Then two cases
are possible.
\\[.1cm]
(1)\quad Suppose $i'<i$; that is, $0=k'_i<k_i$. Then, by Remark~6,
every pair of $t'$ must be of the form $(V'W',U'W')$ for
$|V'|=|U'|=i'+1$ and $|W'|=n-i'-1$. On the other hand, by
Remark~7 every diagonal of $t$ is of the form $(VLW,URW)$ or
$(VRW,ULW)$ for $|V|=|U|=i$ and $|W|=n-i-1$. Since
$|W|<|W'|$, every diagonal of $t$ is not a pair of $t'$.
\\[.1cm]
(2)\quad Suppose $i'=i$; that is, $0<k'_i<k_i$. By Lemma~3, the second
member of every diagonal of $t$ is of one of the following forms:
$VL(LR)^l W$ and $VR(RL)^l W$ for $k_i=2l$, $VRR(LR)^lW$ and $VLL(RL)^lW$ for
$k_i=2l+1$, where $|V|=i-k_i$ and $|W|=n-i-1$. On the other hand, by
Remark~4, no pair of $t'$ can have such second members.

Suppose now for the induction step that $\nu>0$. If $\max(t)=j$ and $\max(t')=j'$,
then $j=j'\geq 2$ and $k_j=k'_j$. Every diagonal $(W_1,W_3)$ of $t$
can be obtained by composing a diagonal $(W_1,W_2)$ of
$h^n_{1,k_1}\ldots h^n_{j-1,k_{j-1}}$ and a diagonal $(W_2,W_3)$ of
$h^n_{j,k_j}$.
By the induction hypothesis,
in $h^n_{1,k'_1}\ldots h^n_{j-1,k'_{j-1}}h^n_{j,0}\ldots h^n_{n-1,0}$
we don't have the pair $(W_1,W_2)$, which is a diagonal of
$h^n_{1,k_1}\ldots h^n_{j-1,k_{j-1}}h^n_{j,0}\ldots h^n_{n-1,0}$,
since the former expression precedes the latter in the lexicographical order.

Suppose $(W_1,W_3)$ is in $t'$. Then for some $W_4$
we need to have $(W_1,W_4)$ in
$h^n_{1,k'_1}\ldots h^n_{j-1,k'_{j-1}}$ and $(W_4,W_3)$ in
$h^n_{j,k_j}$. It is excluded that $W_4$ is $W_2$, as we remarked above.

Suppose $k_j=2l+1$, and suppose $(W_2,W_3)$ is $(V(LR)^l LLW,VRR(LR)^l W)$
for $|W|=n-j-1$. By Lemma~3, this is one possible case. Then, by
Remark~5, we must have that $W_4$ is $W'_4S_qW$ for $S_q$ different
from $L$, since $W_4$ cannot be $W_2$. By Remark~7, we have that
$W_1$ must be of the form
$W'_1LW$. By Remark~6, we cannot have $(W'_1LW,W'_4S_qW)$ in
$h^n_{1,k'_1}\ldots h^n_{j-1,k'_{j-1}}$. So $(W_1,W_3)$ is not in $t'$.
We proceed analogously in other possible cases.
\qed

As a corollary of Lemma~8 and of the Theorem, we obtain that every matrix in \Jn\ has
an entry with 1 (one of its ``diagonals'') where all the matrices preceding it
in the lexicographical order have 0. This implies that the matrices assigned to
the elements of \Jn\ make a linearly independent list, and hence
Brauer's representation
of Temperley-Lieb algebras in matrices, which we have considered here,
is faithful.

\end{document}